\documentclass[12pt]{article}
\usepackage{latexsym, amssymb, amscd, amsmath}
\begin{document}
\title{\textbf{Brill-Noether Loci and the \\ Gonality Stratification of $\mathcal{M}_g$ }}
\author{GAVRIL FARKAS
}
\date{ }
\maketitle 

\newtheorem{thm}{Theorem}
\newtheorem{prop}{Proposition}[section]

\section{Introduction}
\noindent For an irreducible smooth projective complex curve $C$ of genus $g$, the gonality defined as $\mbox{gon}(C)=\mbox{min}\{d\in \mathbb Z_{\geq 1}:\mbox{there exists a }\mathfrak g^1_d\mbox{ on }C\}$ is perhaps the second most natural invariant: it gives an indication of how far $C$ is from being rational, in a way different from what the genus does. For $g\geq 3$ we consider the stratification of the moduli space $\mathcal{M}_g$ of smooth curves of genus $g$ given by gonality:
$$\mathcal{M}_{g,2}^1\subseteq \mathcal{M}_{g,3}^1\subseteq \ldots \subseteq \mathcal{M}^1_{g,k}\subseteq \ldots \subseteq \mathcal{M}_g,$$
where $\mathcal{M}^1_{g,k}:=\{[C]\in \mathcal{M}_g:C\mbox{ has a }\mathfrak g^1_k\}.$ It is well-known that the $k$-gonal locus $\mathcal{M}_{g,k}^1$ is an irreducible variety of dimension $2g+2k-5$ when $k\leq (g+2)/2$; when $k\geq [(g+3)/2]$ one has that $\mathcal{M}^1_{g,k}=\mathcal{M}_g$ (see for instance [AC]). The number
$[(g+3)/2]$ is thus the {\sl generic gonality} for curves of genus $g$.
\newline
\indent For positive integers $g,r$ and $d$, we introduce the {\sl Brill-Noether} locus 
$$\mathcal{M}^r_{g,d}=\{[C]\in \mathcal{M}_g:C\mbox{ carries a }\mathfrak g^r_d\}.$$
The Brill-Noether Theorem (cf. [ACGH]) asserts that when the {\sl Brill-Noether number} $\rho(g,r,d)=g-(r+1)(g-d+r)$ is negative, the general curve of genus $g$
has no $\mathfrak g^r_d$'s, hence in this case the locus $\mathcal{M}_{g,d}^r$ is a proper subvariety of $\mathcal{M}_g$. We study the relative position of the
loci $\mathcal{M}_{g,d}^r$ when $r\geq 3$ and $\rho(g,r,d)<0$ with respect to the gonality stratification of $\mathcal{M}_g$. Typically, we would like to know the gonality of a `general' point $[C]\in \mathcal{M}_{g,d}^r$, or equivalently
the gonality of a `general' smooth curve $C\subseteq \mathbb P^r$ of genus $g$ and degree $d$. Since the geometry of the loci $\mathcal{M}^r_{g,d}$ 
is very messy (existence of many components, some nonreduced and/or not of expected dimension), we will content ourselves with computing $\mbox{gon}(C)$ when $[C]$ is a general point of a
`nice' component of $\mathcal{M}_{g,d}^r$ (i.e. a component which is generically smooth, of the expected dimension and with general point corresponding to a curve with a very
ample $\mathfrak g^r_d$).
\newline
\indent Our main result is the following:
\begin{thm}
Let $g\geq 15$ and $d\geq 14$ be integers with $g$ odd and $d$ even, such that \ $d^2>8g$, $4d<3g+12$, $d^2-8g+8$ is not a square and either $d\leq 18$ or $g<4d-31$. If $$(d',g')\in \{(d,g),(d+1,g+1),
(d+1,g+2),(d+2,g+3)\},$$ then there exists a regular component of the Hilbert scheme $\rm{Hilb}$$_{d',g',3}$ whose general point $[C']$ is a smooth curve such that
$\rm{gon}$$(C')=$$\rm{min}$$(d'-4,[(g'+3)/2])$.
\end{thm}
Here by $\mbox{Hilb}_{d,g,r}$ we denote the Hilbert scheme of curves $C\subseteq \mathbb P^r$ with $p_a(C)=g$ and $\mbox{deg}(C)=d$. A component of $\mbox{Hilb}_{d,g,r}$ is said to be {\sl regular} if its general point corresponds to a smooth irreducible curve $C\subseteq \mathbb P^r$ such that the
normal bundle $N_{C/\mathbb P^r}$ satisfies $H^1(C,N_{C/\mathbb P^r})=0$. By standard deformation theory (cf. [Mod] or [Se]), a regular component of $\mbox{Hilb}_{d,g,r}$ is generically smooth of the expected dimension $\chi(C,N_{C/\mathbb P^r})=(r+1)d-(r-3)(g-1)$. Note that for $r=3$ the expected dimension of the Hilbert scheme is just $4d$. We refer to Section 4 for a natural extension of Theorem 1 for curves in higher dimensional projective spaces.
\newline
\indent As for the numerical conditions entering Theorem 1, we note that the inequality $d^2>8g$ ensures the existence of smooth curves $C\subseteq \mathbb P^3$ with $g(C)=g$ and $\mbox{deg}(C)=d$ (see Section 2), $4d<3g+12\Leftrightarrow \rho(g,3,d)<0$ is just the condition that $\mathcal{M}_{g,d}^3$ is a proper subvariety of $\mathcal{M}_g$, while the remaining requirements are mild technical conditions.
\newline
\indent A remarkable application of Theorem 1 is a new proof of our result (cf. [Fa]):
\begin{thm}
The Kodaira dimension of the moduli space of curves of genus $23$ is $\geq 2$.
\end{thm}
We recall that for $g\geq 24$ Harris, Mumford and Eisenbud proved (cf. [HM],[EH]) that $\mathcal{M}_g$ is of general type whereas for $g\leq 16, g\neq 14$ we have that $\kappa(\mathcal{M}_g)=-\infty.$ The famous Slope Conjecture of
Harris and Morrison predicts that $\mathcal{M}_g$ is uniruled for all $g\leq 22$
(see [Mod]). Therefore the moduli space $\mathcal{M}_{23}$ appears as an intriguing transition case between two extremes: uniruledness and being of general type.\vskip 5pt
\indent To put our Theorem 1 into perspective, let us note that for $r=2$ we have the following result of M. Coppens (cf. [Co]): let $\nu:C\rightarrow \Gamma$ be the normalization of a general, irreducible plane curve of degree $d$ with $\delta =g-{d-1\choose 2}$ nodes. Assume that $0<\delta <(d^2-7d+18)/2$. Then $\rm{gon}$$(C)=d-2.$
\newline
\indent This theorem says that there are no $\mathfrak g^1_{d-3}$'s
on $C$. On the other hand a $\mathfrak g^1_{d-2}$ is given by the lines through a node of $\Gamma $. The condition $\delta <(d^2-7d+18)/2$ from the statement is equivalent with $\rho(g,1,d-3)<0$. This is the range in which the problem is non-trivial: if $\rho(g,1,d-3)\geq 0$, the Brill-Noether Theorem provides $\mathfrak g^1_{d-3}$'s on $C$. 
\vskip 5pt
\indent For $r\geq 3$ we might hope for a similar result. Let $C\subseteq \mathbb P^r$ be a suitably general smooth curve of genus $g$ and degree $d$, with $\rho(g,r,d)<0$. We can always assume that $d\leq g-1$ (by duality $\mathfrak g^r_d\mapsto |K_C-\mathfrak g^r_d|$ we can always land in this range). One can expect that a $\mathfrak g^1_k$ computing
$\mbox{gon}(C)$ is of the form $\mathfrak g^r_d(-D)=\{ E-D:E\in \mathfrak g^r_d, E\geq D\}$ for some effective divisor $D$ on $C$. Since the expected dimension of the variety of $e$-secant $(r-2)$-plane divisors
$$V_e^{r-1}(\mathfrak g^r_d):=\{ D\in C_e:\mbox{ dim } \mathfrak g^r_d(-D)\geq 1\}$$
is $2r-2-e$ (cf. [ACGH]), we may ask whether $C$ has finitely many $(2r-2)$-secant $(r-2)$-planes (and no $(2r-1)$-secant $(r-2)$-planes at all). This is known to be true for curves with general moduli, that is, when $\rho(g,r,d)\geq 0$ (cf. [Hir]): for instance a smooth curve $C\subseteq \mathbb P^3$ with
general moduli has only finitely many $4$-secant lines and no $5$-secant lines. No such principle appears to be known for curves
with special moduli.\vskip 4pt
\noindent \textbf{Definition: } We call the number $\mbox{min}(d-2r+2,[(g+3)/2])$ the {\sl expected gonality } of a smooth nondegenerate curve $C\subseteq \mathbb P^r$ of degree $d$ and genus $g$.
\newline 
\indent One can approach such problems from a different angle: find recipes to compute the gonality of various classes of curves $C\subseteq \mathbb P^r$. Our
knowledge in this respect is very scant: we know how to compute the gonality of extremal curves $C\subseteq \mathbb P^r$ (that is, curves attaining the Castelnuovo bound, see [ACGH]) and the gonality of complete intersections in $\mathbb P^3$ (cf. [Ba]): If $C\subseteq \mathbb P^3$ is a smooth complete intersection of type $(a,b)$ then $\mbox{gon}(C)=ab-l$, where $l$ is the degree of a maximal
linear divisor on $C$. Hence an effective divisor $D\subseteq C$ computing
$\mbox{gon}(C)$ is residual to a linear divisor of degree $l$ in a
plane section of $C$. 
\newline
\noindent{\textbf{Acknowledgments:}} This paper is part of my thesis written at the Universiteit van Amsterdam. The help of my advisor Gerard van der Geer, and of Joe Harris, is gratefully acknowledged.

\section{Linear systems on $K3$ surfaces in $\mathbb P^r$}

\noindent We will construct smooth curves $C\subseteq \mathbb P^r$ having the expected gonality starting with sections of smooth $K3$ surfaces. We recall a few basic facts about linear systems on $K3$ surfaces
(cf. [SD]).\newline
\indent Let $S$ be a smooth $K3$ surface. For an effective
divisor $D\subseteq S$, we have $h^1(S,D)=h^0(D, \mathcal{O}_D)-1$. If $C\subseteq S$ is an irreducible curve then $H^1(S,C)=0$, and by Riemann-Roch we have that
$\mbox{dim}|C|=1+C^2/2=p_{a}(C).$
In particular $C^2\geq -2$ for every irreducible curve $C$. Moreover we have equivalences \newline \indent
$C^2=-2 \Longleftrightarrow \mbox{dim}|C|=0 \Longleftrightarrow C \mbox{ is a smooth rational curve }\mbox{ and }$\vskip 4pt \indent
$C^2=0 \Longleftrightarrow  \mbox{dim}|C|=1 \Longleftrightarrow p_a(C)=1.$\vskip 5pt
\indent
For a $K3$ surface one also has a `strong Bertini' Theorem (cf. [SD]):
\begin{prop}
Let $\mathcal{L}$ be a line bundle on a $K3$ surface $S$ such that $|\mathcal{L}|\neq \emptyset$. Then $|\mathcal{L}|$
has no base points outside its fixed components. Moreover, if $\rm{bs }|\mathcal{L}|=\emptyset$ then either
\begin{itemize}
\item $\mathcal{L}^2> 0,\mbox{ }h^1(S,\mathcal{L})=0$ and the general member of $|\mathcal{L}|$ is a smooth, irreducible curve of genus $\mathcal{L}^2/2+1$,
or
\item $\mathcal{L}^2=0$ and $\mathcal{L}=\mathcal{O}_S(kE)$, where $k\in \mathbb Z_{\geq 1}$,  $E\subseteq S$ is an irreducible curve with $p_a(E)=1$. 
We have that $h^0(S,\mathcal{L})=k+1, \mbox{ } h^1(S,\mathcal{L})=k-1$
and all divisors in $|\mathcal{L}|$ are of the form $E_1+\cdots +E_k$ with $E_i\sim E$.
\end{itemize}
\end{prop}
\indent We are interested in space curves sitting on $K3$ surfaces and the
starting point is Mori's Theorem (cf. [Mo]): if $d>0$, $g\geq 0$, there is a
smooth curve $C\subseteq \mathbb P^3$ of degree $d$ and genus $g$, lying on a 
smooth quartic surface $S$, if and only if (1) $g=d^2/8+1$, or (2) $g<d^2/8$ and $(d,g)\neq (5,3)$. Moreover, we can choose $S$ such that $\mbox{Pic}(S)=\mathbb Z H=\mathbb Z (4/d)C$ in case (1) and such that $\mbox{Pic}(S)=
\mathbb Z H\oplus \mathbb Z C$, with $H^2=4, C^2=2g-2$ and $H\cdot C=d$, in case (2). In each case $H$ denotes a plane section of $S$. Note that
from the Hodge Index Theorem one has the necessary condition 
$(C\cdot H)^2-H^2 C^2=d^2-8(g-1)\geq 0.$\newline
\indent Mori's result has been extended by Rathmann to curves in higher dimensional projective spaces (cf. [Ra], see also [Kn]): For integers $d>0, g>0$ and $r\geq 3$ such that $d^2\geq 4g(r-1)+(r-1)^2$, there exists a smooth $K3$ surface $S\subseteq \mathbb P^r$ of degree $2r-2$ and a smooth curve $C\subseteq S$ of genus $g$ and degree $d$ such that $\mbox{Pic}(S)=\mathbb Z H\oplus \mathbb Z C$, where $H$ is a hyperplane section of $S$.
\newline
\indent We will repeatedly use the following simple observation:
\begin{prop}
Let $S\subseteq \mathbb P^r$ be a smooth $K3$ surface of degree $2r-2$ with a smooth curve $C\subseteq S$ such that $\rm{Pic}$$(S)=\mathbb Z H\oplus \mathbb Z C$ and assume that $S$ has no $(-2)$ curves. A divisor class $D$ on $S$ is effective if and only if $D^2\geq 0$ and $D\cdot H>2$.
\end{prop}
\vskip 5pt
\textbf{Remark:} If $S\subseteq \mathbb P^r$ is a smooth $K3$ surface of degree $2r-2$ with Picard number $2$ as above, $S$ has no $(-2)$ curves when the equation
\begin{equation}
(r-1)m^2+mnd+(g-1)n^2=-1
\end{equation}
has no solutions $m,n\in \mathbb Z$. This is the case for instance when $d$ is even and $g$ and $r$ are odd. Furthermore, a necessary condition for $S$ to have genus $1$ curves is that $d^2-4(g-1)(r-1)$ is a square.

\section{Brill-Noether special linear series on curves on $K3$ surfaces}
\noindent The first important result in the study of special linear series on curves lying on $K3$ surfaces
was Lazarsfeld's proof of the Brill-Noether-Petri Theorem (cf. [Laz]). He noticed that there is no Brill-Noether type obstruction to embed a curve in a $K3$
surface: if $C_0\subseteq S$ is a smooth curve of genus $g\geq 2$ on a $K3$
surface such that $\mbox{Pic}(S)=\mathbb ZC_0$, then the general curve $C\in |C_0|$ satisfies the Brill-Noether-Petri Theorem, that is, for any line bundle
$A$ on $C$, the Petri map $\mu_0(C,A):H^0(C,A)\otimes H^0(C,K_C\otimes A^{\vee})\rightarrow H^0(C,K_C)$ is injective. We mention that Petri's Theorem 
implies (trivially) the Brill-Noether Theorem.  
\newline
\indent The general philosophy when studying linear series on a $K3$-section
$C\subseteq S$ of genus $g\geq 2$, is that the type of a Brill-Noether special
$\mathfrak g^r_d$ often does not depend on $C$ but only on its linear equivalence class in $S$, i.e. a $\mathfrak g^r_d$ on $C$ with $\rho(g,r,d)<0$
is expected to propagate to all smooth curves $C'\in |C|$. This expectation, in such generality, is perhaps a bit too optimistic, but it was proved to be true
for the Clifford index of a curve (see [GL]): for $C\subseteq S$ a smooth $K3$-section of genus $g\geq 2$, one has that $\mbox{Cliff}(C')=\mbox{Cliff}(C)$
for every smooth curve $C'\in |C|$. Furthermore, if $\mbox{Cliff}(C)<[(g-1)/2]$
(the generic value of the Clifford index), then there exists a line bundle
$\mathcal{L}$ on $S$ such that for all smooth $C'\in |C|$ the restriction $\mathcal{L}_{|C'}$ computes $\mbox{Cliff}(C')$. Recall that the {\sl Clifford index} of a curve $C$ of genus $g$ is defined as 
$$\mbox{Cliff}(C):=\mbox{min}\{\mbox{Cliff}(D): D\in \mbox{Div}(C),h^0(D)\geq 2, h^1(D)\geq 2\},$$
where for an effective divisor $D$ on $C$, we have $\mbox{Cliff}(D)=\mbox{deg}(D)-2(h^0(D)-1)$.
Note that in the definition of $\mbox{Cliff}(C)$ the condition $h^1(D)\geq 2$ can be replaced with $\mbox{deg}(D)\leq g-1$.
Another invariant of a curve is the {\sl Clifford dimension} of $C$ defined as $$\mbox{Cliff-dim}(C):=\mbox{min}\{ r\geq 1:\exists \mathfrak g^r_d \mbox{ on } C\mbox{ with }d\leq g-1,\mbox{ such that } d-2r=\mbox{Cliff}(C)\}.$$ 
Curves with Clifford dimension $\geq 2$ are rare: smooth plane curves are precisely the curves of Clifford dimension $2$, while curves of Clifford dimension $3$ occur only in genus $10$ as complete
intersections of two cubic surfaces in $\mathbb P^3$.
\newline
\indent Harris and Mumford during their work in [HM] conjectured that the gonality of a $K3$-section should stay constant in a linear system: if $C\subseteq S$ carries an exceptional $\mathfrak g^1_d$ then every smooth $C'\in |C|$ carries an equally
exceptional $\mathfrak g^1_d$. This conjecture was later disproved by Donagi and Morrison (cf. [DMo]). They came up with the following counterexample: let $\pi :S\rightarrow \mathbb P^2$ be a
$K3$ surface, double cover of $\mathbb P^2$ branched along a smooth sextic
and let $\mathcal{L}=\pi^*\mathcal{O}_{\mathbb P^2}(3)$. The genus of a smooth $C\in |\mathcal{L}|$ is $10$. The general $C\in |\mathcal{L}|$ carries a very ample $\mathfrak g^2_6$, hence $\mbox{gon}(C)=5$. On the other hand, any curve
in the codimension $1$ linear system $|\pi^*H^0(\mathbb P^2,\mathcal{O}_{\mathbb P^2}(3))|$ is bielliptic, therefore has gonality $4$.
Under reasonable assumptions this turns out to be the only counterexample to the Harris-Mumford conjecture. Ciliberto and Pareschi
proved that if $C\subseteq S$ is such that $|C|$ is base-point-free
and ample, then either $\mbox{gon}(C')=\mbox{gon}(C)$ for all smooth $C'\in |C|$, or $(S,C)$ are as in the previous counterexample (cf. [CilP]).
\newline
\indent Although $\mbox{gon}(C)$ can drop as $C$ varies in a linear system, base-point-free $\mathfrak g^1_d$'s on $K3$-sections do propagate:
\begin{prop}[Donagi-Morrison]
Let $S$ be a $K3$ surface, $C\subseteq S$ a smooth, nonhyperelliptic curve and
$|Z|$ a complete, base-point-free $\mathfrak{g}^1_d$ on $C$ such that $\rho(g,1,d)<0$. Then there is an effective divisor $D\subseteq S$ such that:
\begin{itemize}
\item $h^0(S,D)\geq 2,\mbox{ h}^0(S,C-D)\geq 2,\rm{ deg}_{C}(D_{|C})\leq g-1.$
\item $\rm{Cliff}$$(C',D_{|C'})\leq \rm{Cliff}$$(C,Z)$, for any smooth $C'\in |C|$.
\item There is $Z_0\in |Z|$, consisting of distinct points such that $Z_0\subseteq D\cap C$.
\end{itemize}
\end{prop}
\noindent Throughout this paper, for a smooth curve $C$ we denote, as usual, by $W^r_d(C)$ the scheme whose points are line bundles $A\in \mbox{Pic}^d(C)$ with $h^0(C,A)\geq r+1$, and by $G^r_d(C)$ the scheme parametrizing $\mathfrak g^r_d$'s on $C$.

\section{The gonality of curves in $\mathbb P^r$}

For a wide range of $d,g$ and $r$ we construct curves $C\subseteq \mathbb P^r$ of degree $d$ and genus $g$ having the expected gonality. We start with a case when we can realize our curves as sections of $K3$ surfaces.
\begin{thm}
Let $r\geq 3, d\geq r^2+r$ and $g\geq 0$ be integers such that $\rho(g,r,d)<0$ and with $d^2>4(r-1)(g+r-2)$ when $r\geq 4$ while $d^2>8g$ when $r=3$. Let us assume moreover that $0$ and $-1$ are not represented by the quadratic form 
$$Q(m,n)=(r-1)m^2+mnd+(g-1)n^2,\mbox{ }\mbox{ }m,n\in \mathbb Z.$$ 
Then there exists a smooth curve $C\subseteq \mathbb P^r$ of degree $d$ and genus $g$ such that $\rm{gon}$$(C)=\rm{min}$$(d-2r+2,[(g+3)/2])$. If $\rm{gon}$$(C)=d-2r+2<[(g+3)/2]$ then $\rm{dim}$$\mbox{ }W^1_{d-2r+2}(C)=0$ and every $\mathfrak g^1_{d-2r+2}$ is given by the hyperplanes through a $(2r-2)$-secant 
$(r-2)$-plane.
\end{thm}

\noindent {\sl Proof: } By Rathmann's Theorem there exists a smooth $K3$ surface $S\subseteq \mathbb P^r$ with $\mbox{deg}(S)=2r-2$ and $C\subseteq S$ a smooth curve of degree $d$ and genus $g$ such that $\mbox{Pic}(S)=\mathbb Z H\oplus \mathbb Z
C$, where $H$ is a hyperplane section. The conditions $d,g$ and $r$ are subject to, ensure that $S$ does not contain $(-2)$ curves or genus $1$ curves.
\newline
\indent We prove first that $\mbox{Cliff-dim}(C)=1$. It suffices to show that $C\subseteq S$ is an ample divisor, because then by using Prop.3.3
from [CilP] we obtain that either $\mbox{Cliff-dim}(C)=1$ or $C$ is a smooth plane
sextic, $g=10$ and $(S,C)$ are as in Donagi-Morrison's example (then $\mbox{Cliff-dim}(C)=2$). The latter case obviously does not happen. 
\newline
\indent We prove that $C\cdot D>0$ for any effective
divisor $D\subseteq S$. Let $D\sim mH+nC$, with $m,n\in \mathbb Z$, such a
divisor. Then $D^2=(2r-2)m^2+2mnd+n^2(2g-2)\geq 0$ and $D\cdot H=(2r-2)m+dn>2.$ The case $m\leq 0,n\leq 0$ is impossible, while the case $m\geq 0,n\geq 0$ is trivial.
Let us assume $m>0,n<0$. Then $D\cdot C=md+n(2g-2)>-n\bigl(d^2/(2r-2)-2g+2\bigr )+d/(r-1)>0,$ because
$d^2/(2r-2)>2g$. In the remaining case $m<0,n>0$ we have that $nD\cdot C\geq -mD\cdot H>0$,
so $C$ is ample by Nakai-Moishezon.
\newline 
\indent Our assumptions imply that $d\leq g-1$, so $\mathcal{O}_C(1)$ is among the line bundles from which $\mbox{Cliff}(C)$ is computed. We get thus the following estimate on the gonality of $C$:

$$\mbox{gon}(C)=\mbox{Cliff}(C)+2\leq \mbox{Cliff}(C,H_{|C})+2=d-2r+2,$$
which yields $\mbox{gon}(C)\leq \mbox{min}(d-2r+2,[(g+3)/2]).$ 
\newline
\indent For the rest of the proof let us assume that $\mbox{gon}(C)<[(g+3)/2]$. We will then show that $\mbox{gon}(C)=d-2r+2.$ Let $|Z|$ be a complete, base point free pencil computing $\mbox{gon}(C).$ By applying Prop.3.1, there exists an effective divisor
$D\subseteq S$ satisfying
$$h^0(S,D)\geq 2, h^0(S,C-D)\geq 2,
\mbox{deg}(D_{|C})\leq g-1, \mbox{ gon}(C)=\mbox{Cliff}(D_{|C})+2\mbox{ and }Z\subseteq D\cap C.$$ 
We consider the exact cohomology sequence:
$$0\rightarrow H^0(S,D-C)\rightarrow H^0(S,D)\rightarrow H^0(C,D_{|C})\rightarrow H^1(S,D-C).$$
Since $C-D$ is effective and $\nsim 0$, one sees that $D-C$ cannot be effective, so $H^0(S,D-C)=0$. The surface $S$ does not contain $(-2)$ curves, so
$|C-D|$ has no fixed components; the equation $(C-D)^2=0$ has no solutions,
therefore $(C-D)^2>0$ and the general element of $|C-D|$ is smooth and irreducible. Then it follows that $H^1(S,D-C)=H^1(S,C-D)^{\vee}=0.$
 Thus $H^0(S,D)=H^0(C,D_{|C})$ and
$$\mbox{gon}(C)=2+\mbox{Cliff}(D_{|C})=2+D\cdot C-2\mbox{ dim}|D|=D\cdot C-D^2.$$
We consider the following family of effective divisors 

$$\mathcal{A}:=\{D\in \mbox{Div}(S):h^0(S,D)\geq 2,h^0(S,C-D)\geq 2,\mbox{ } C\cdot D\leq g-1\}.$$
Since we already know that $d-2r+2\geq \mbox{gon}(C)\geq \alpha $, where
$\alpha=\mbox{min}\{D\cdot C-D^2:D\in \mathcal{A}\}$, we are done if we prove that $\alpha\geq d-2r+2$. Take $D\in \mathcal{A}$ such that $D\sim mH+nC$, $m,n\in \mathbb Z$. The conditions $D^2>0, D\cdot C\leq g-1$ and $2<D\cdot H<d-2$ (use Prop.2.2
for the last inequality) can be rewritten as 
$$(r-1)m^2+mnd+n^2(g-1)>0\mbox{ (i), } 2<(2r-2)m+nd<d-2\mbox{ (ii), }md+(2n-1)(g-1)\leq
0 \mbox{ (iii)}.$$
We have to prove that for any $D\in \mathcal{A}$ the following inequality
holds
$$f(m,n)=D\cdot C-D^2=-(2r-2)m^2+m(d-2nd)+(n-n^2)(2g-2) \geq f(1,0)=d-2r+2.$$ 
We solve this standard calculus problem. Denote by $$a:=\frac{d+\sqrt{d^2-4(r-1)(g-1)\ }}{2r-2}\mbox{ }\mbox{ and }\mbox{ }b:=\frac{d-\sqrt{d^2-4(r-1)(g-1)\ }}{2r-2}\mbox{ }.$$ We dispose first of the case $n<0$. Assuming $n<0$, from (i) we
have that either $m<-bn$ or $m>-an$. If $m<-bn$ from (ii) we obtain that 
$2<n(d-(2r-2)b)<0$, because $n<0$ and $d-(2r-2)b=\sqrt{d^2-4(r-1)(g-1)}>0$, so we have reached a contradiction. \newline
\indent We assume now that $n<0$ and $m>-an$.  From (iii) we get that $m\leq (g-1)(1-2n)/d$.
If $-an>(g-1)(1-2n)/d$ we are done because there are no $m,n\in \mathbb Z$ satisfying (i), (ii) and (iii), while in the other case for any $D\in \mathcal{A}$ with $D\sim mH+nC$, one has the inequality
$$f(m,n)>f(-an,n)=\frac{(d^2-4(r-1)(g-1))+d\sqrt{d^2-4(r-1)(g-1)}}{2r-2}(-n).$$
\indent When $r\geq 4$ since we assume that $\sqrt{d^2-4(r-1)(g-1)}\geq 2r-2$,
it immediately follows that $f(m,n)\geq d>d-2r+2.$ In the case $r=3$ when we only have the weaker assumption $d^2>8g$, we still get that $f(-an,n)>d-4$ 
unless $n=-1$ and $d^2-8g<8$. In this last situation we 
obtain $m\geq (d+4)/4$ so $f(m,-1)\geq f((d+4)/4,-1)>d-4$.
\newline
\indent The case $n>0$ can be treated rather similarly. From (i) we get that either $m<-an$ or $m>-bn$. The first case can be dismissed immediately. When
$m>-bn$ we use that for any $D\in \mathcal{A}$ with  $D\sim mH+nC$,
$$f(m,n)\geq \mbox{min}\bigl\{f(-(g-1)(2n-1)/d,n),\mbox{max}\{f(-bn,n),f\bigl((2-nd)/(2r-2),n\bigr)\}\bigr \}.$$
Elementary manipulations give that 
$$f(-(g-1)(2n-1)/d,n)=(g-1)/2\ [(2n-1)^2(d^2-4(r-1)(g-1))/{d^2}+1]\geq d-2r+2$$
(use only that $d\leq g-1$ and $d^2>4(r-1)g$, so we cover both cases $r=3$ and $r\geq 4$ at once). Note that in the case $n>0$ we have equality if and only if $n=1,m=-1$ and $d=g-1$.
\newline
\indent Moreover $f(-bn,n)=n(2g-2-bd)\geq 2g-2-bd$
and $2g-2-bd> d-2r+2\ \Leftrightarrow 2r-2< \sqrt{d^2-4(r-1)(g-1)}< d-2r+2$. When this
does not happen we proceed as follows: if $\sqrt{d^2-4(r-1)(g-1)}\geq d-2r+2$ then if $n=1$
we have that $m>-b\geq -1$, that is $m\geq 0$, but this contradicts (ii). When $n\geq 2$, we have $f((2-nd)/(2r-2),n)=[(d^2-4(r-1)(g-1))(n^2-n)+(2d-4)]/(2r-2)> d-2r+2.$ Finally,
the remaining possibility $2r-2\geq \sqrt{d^2-4(r-1)(g-1)}$ does not occur when $r\geq 4$ while in the case $r=3$ we either have $f(-bn,n)>d-4$ or else $n=1$ and then $m>(-d+4)/4$ hence $f(m,1)>f((-d+4)/4,1)=d-4.$ 
\newline
\indent All this leaves us with the case $n=0$, when $f(m,0)=-(2r-2)m^2+md$. Clearly
$f(m,0)\geq f(1,0)$ for all $m$ complying with (i),(ii) and (iii). 
\newline
\indent Thus we proved that $\mbox{gon}(C)=d-2r+2$. We have equality $D\cdot C-D^2=d-2r+2$ where $D\in \mathcal{A}$, if and only if $D=H$ or in the case $d=g-1$ also when $D=C-H$. The latter possibility can be ruled out since $d=g-1$ is not compatible with the assumptions $d\geq r^2+r$ and $d-2r+2<[(g+3)/2]$. Therefore we can always assume that the divisor on $S$ cutting a $\mathfrak g^1_{d-2r+2}$ on $C$ is the hyperplane section of $S$. Since $Z\subseteq H\cap C$, if we denote by $\Delta $ the residual divisor 
of $Z$ in $H\cap C$, we have that $h^0(C,H_{|C}-\Delta)=2$, so $\Delta$
spans a $\mathbb P^{r-2}$ hence $|Z|$ is given by the hyperplanes through the $(2r-2)$-secant $(r-2)$-plane $\langle \Delta \rangle$. This shows that every pencil computing $\mbox{gon}(C)$ is given by the hyperplanes through a $(2r-2)$-secant $(r-2)$-plane.
\newline
\indent There are a few ways to see that $C$ has only finitely many $(2r-2)$-secant $(r-2)$-planes. The shortest is to invoke Theorem 3.1 from [CilP]: since $\mbox{gon}(C')=d-2r+2$ is constant as $C'$ varies in $|C|$, for the general smooth curve $C'\in |C|$ one has $\mbox{dim }W^1_{d-2r+2}(C')=0$. \hfill $\Box$
\vskip 5pt
\noindent \textbf{Remarks: 1.} Keeping the assumptions and the notations of Theorem 3 we note that when $d-2r+2<[(g+3)/2]$ the linear system $|C|$ is $(d-2r-1)$-very ample, i.e. for any $0$-dimensional subscheme $Z\subseteq S$ of length $\leq d-2r$ the map $H^0(S,C)\rightarrow H^0(S,C\otimes \mathcal{O}_Z)$\  is surjective. Indeed, by applying Theorem 2.1 from [BS] if $|C|$ is not $(d-2r-1)$-very ample, there exists an effective divisor $D$
on $S$ such that $C-2D$ is $\mathbb Q$-effective and 
$$C\cdot D-(d-2r)\leq D^2\leq C\cdot D/2<d-2r,$$
hence $C\cdot D-D^2\leq d-2r$. On the other hand clearly $D\in \mathcal{A}$,
thus $C\cdot D-D^2\geq d-2r+2$, a contradiction.
\newline
\noindent \textbf{2. }
One can find quartic surfaces $S\subseteq \mathbb P^3$ containing a smooth curve $C$ of degree $d$ and genus $g$ in the case $g=d^2/8+1$ (which is outside the range Theorem 3 deals with). Then $d=4m,g=2m^2+1$ with $m\geq 1$ and $C$ is a complete intersection of type $(4,m)$. For such a curve, $\mbox{gon}(C)=d-l$, where $l$ is the degree of a maximal linear divisor on $C$ (cf. [Ba]). If $S$ is sufficiently general so that it contains no lines, by Bezout, $C$ cannot have $5$-secant lines so $\mbox{gon}(C)=d-4$ in this case too.\newline
\vskip 8pt
When $r=3$ we want to find out when the curves constructed in Theorem 3 correspond to smooth
points of $\mbox{Hilb}_{d,g,3}$. We have the following:
\begin{prop}
Let $C\subseteq S\subseteq \mathbb P^3$ be a smooth curve sitting on a quartic surface such that $\rm {Pic}$$(S)=\mathbb Z H\oplus \mathbb Z C$ with $H$ being a 
plane section and assume furthermore that $S$ contains no $(-2)$ curves.
Then $H^1(C,N_{C/\mathbb P^3})=0$ if and only if $ d\leq 18 \mbox{ or }g<4d-31$.
\end{prop}
{\sl Proof:} We use the exact sequence
\begin{equation}
0\longrightarrow N_{C/S}\longrightarrow N_{C/\mathbb P^3}\longrightarrow N_{S/\mathbb P^3}\otimes \mathcal{O}_C\longrightarrow 0,
\end{equation}
where $N_{S/\mathbb P^3}\otimes \mathcal{O}_C=\mathcal{O}_C(4)$ and $N_{C/S}=K_C$. We claim that there is an isomorphism $H^1(C,N_{C/\mathbb P^3})=H^1(C,\mathcal{O}_C(4))$. Suppose this is not the case. Then the injective map
$H^1(C,K_C)\rightarrow H^1(C,N_{C/\mathbb P^3})$ provides a section $\sigma \in H^0(N_{C/\mathbb P^3}^{\vee}\otimes K_C)$ which yields a splitting of the dual of the exact sequence (2), hence (2) is split as well. Using a result from [GH, p.252] we obtain that $C$ is a complete intersection with $S$. This is clearly a contradiction. Therefore one has $H^1(C,N_{C/\mathbb P^3})=H^1(C,\mathcal{O}_C(4))$.
\newline    
\indent We have isomorphisms $H^1(C,4H_{|C})=H^2(S,4H-C)=H^0(S,C-4H)^{\vee}$. According to Prop.2.2 the divisor $C-4H$ is effective if and only if $(C-4H)^2\geq 0$ and $(C-4H)\cdot H>2$, from which the conclusion follows. \hfill $\Box$
\vskip6pt
We need to determine the gonality of nodal curves not of compact type and which 
consist of two components meeting at a number of points. We have the following result:
\begin{prop}
Let $C=C_1\cup_{\Delta}C_2$ be a quasi-transversal union of two smooth curves
$C_1$ and $C_2$ meeting at a finite set $\Delta$. Denote by $g_1=g(C_1), g_2=g(C_2), \delta=\rm{card}(\Delta)$. Let us assume that $C_1$ has only finitely many pencils $\mathfrak g^1_d$, where $\delta \leq d$ and that the points of $\Delta$ do not occur in the same fibre of one of these pencils.
Then $\rm{gon}$$(C)\geq d+1$. Moreover if $\rm{gon}$$(C)=d+1$ then either
(1) $C_2$ is rational and there is a degree $d$ map $f_1:C_1\rightarrow \mathbb P^1$ and a degree $1$ map $f_2:C_2\rightarrow \mathbb P^1$ such that $f_{1\ |\Delta}=f_{2\ |\Delta}$, or (2) there is a $\mathfrak g^1_{d+1}$ on $C_1$ containing $\Delta$ in a fibre.
\end{prop}
{\sl Proof:} Let us assume that $C$ is $k$-gonal, that is, a limit of smooth $k$-gonal curves. If $g=g_1+g_2+\delta-1$, we consider the space $\overline{\mathcal{H}}_{g,k}$ of Harris-Mumford admissible coverings of degree $k$ and we denote by $\pi:\overline{\mathcal{H}}_{g,k}\rightarrow \overline{\mathcal{M}}_g$ the proper map sending a covering to the stable model of its domain (see [HM]). Since $[C]\in \overline{\mathcal{M}}_{g,k}^1=\mbox{Im}(\pi)$,
it follows that there exists a semistable curve $C'$ whose stable model is $C$ and a degree $k$ admissible covering $f:C'\rightarrow Y$, where $Y$ is a semistable curve of arithmetic genus $0$. We thus have that $f^{-1}(Y_{\tiny{sing}})=C'_{\tiny{sing}}$ and if $p\in C'_1\cap C'_2$ with $C'_1$ and $C'_2$ components of $C'$, then $f(C'_1)$ and $f(C'_2)$ are distinct components of $Y$ and the ramification indices at the point $p$ of the restrictions $f_{|C'_1}$ and $f_{|C'_2}$ are the same.
\newline
\indent We have that $C'=C_1\cup C_2 \cup R_1\cup \ldots \cup R_{\delta}$, where for $1\leq i\leq \delta$ the curve $R_i$ is a (possibly empty) destabilizing chain of $\mathbb P^1$'s inserted at the nodes of $C$. Let us denote $\{p_i\}=C_1\cap R_i$ and $\{q_i\}=C_2\cap R_i$; if $R_i=\emptyset$ then we take $p_i=q_i\in \Delta \subseteq C$.
\newline
\indent We first show that $k\geq d+1$. Suppose $k\leq d$. Since $C_1$ has no $\mathfrak g^1_{d-1}$'s it follows that $k=d$ and that $f^{-1}f(C_1)=C_1$.
If there were distinct points $p_i$ and $p_j$ such that $f(p_i)\neq f(p_j)$, then $f(R_i)\neq f(R_j)$ and the image curve $Y$ would no longer have genus $0$.
Therefore $f(p_i)=f(p_j)$ for all $i,j\in \{1,\ldots ,\delta\}$, that is $\Delta$ appears in the fibre of a $\mathfrak g^1_d$ on $C_1$, a contradiction.
\newline
\indent Assume now that $k=d+1$. Then either $\mbox{deg}(f_{|C_1})=d$ or $\mbox{deg}(f_{|C_1})=d+1.$ If $\mbox{deg}(f_{|C_1})=d+1$, then again $f^{-1}f(C_1)=C_1$ and by the same reasoning $f$ maps all the $p_i$'s to the same point and this yields case (2) from the statement of the Proposition. If
$\mbox{deg}(f_{|C_1})=d$ then $f^{-1}f(C_1)=C_1\cup D$, where $D$ is a smooth rational curve mapped isomorphically to its image via $f$. If $D=C_2$ then the condition that the dual graph of $Y$ is a tree implies that $f(p_i)=f(q_i)$ for all $i$
and this yields case (1) from the statement. Finally, if $D\neq C_2$ then $f(C_1)\neq f(C_2)$. We know that there are $1\leq i<j\leq \delta$ such that 
$f(p_i)\neq f(p_j)$. The image $f(C_2)$ belongs to a chain $R$ of $\mathbb P^1$'s such that either $R\cap f(C_1)=\{f(p_i)\}$ or $R\cap f(C_1)=\{f(p_j)\}$.
In the former case $f(p)=f(p_i)$ for all $p\in \Delta-\{p_j\}$ while in the latter case $f(p)=f(p_j)$ for all $p\in \Delta-\{p_i\}.$ In each case by adding a base point we obtain a $\mathfrak g^1_{d+1}$ on $C_1$ containing $\Delta$ in a fibre.\hfill $\Box$  

\vskip6pt
Theorem 3 provides curves $C\subseteq \mathbb P^3$ of expected gonality when $d$ is even and $g$ is odd (equation (1) has no solutions in this case). Naturally, we would like to have such curves when $d$ and $g$ have other parities as well. We will achieve this by attaching to a `good' curve of expected gonality either a $2$ or $3$-secant line or a $4$-secant conic.\vskip 5pt
\noindent {\bf Theorem 1} {\sl Let $g\geq 15$ and $d\geq 14$ be integers with $g$ odd and $d$ even, such that $d^2>8g, 4d<3g+12$, $d^2-8g+8$ is not a square and either $d\leq 18$ or $g<4d-31$. If 
$$(d',g')\in \{(d,g),(d+1,g+1),
(d+1,g+2),(d+2,g+3)\},$$ then there exists a regular component of $\rm{Hilb}_{d',g',3}$ with general point $[C']$ a smooth curve such that
$\rm{gon}$$(C')=\rm{min}(d'-4,[(g'+3)/2])$}.
\vskip 8pt
\noindent {\sl Proof:} For $d$ and $g$ as in the statement we know by Theorem 3 and by Prop.4.1 that there exists a smooth nondegenerate curve $C\subseteq \mathbb P^3$ of degree $d$ and genus $g$, with $\mbox{gon}(C)=\mbox{min}(d-4,[(g+3)/2])$ and $H^1(C,N_{C/\mathbb P^3})=0$. We can also assume that $C$ sits on a smooth quartic surface $S$ and $\mbox{Pic}(S)=\mathbb Z H\oplus \mathbb Z C$. Moreover, in the case $d-4<[(g+3)/2]$ the curve $C$ has only finitely many $\mathfrak g^1_{d-4}$'s, all given by planes through a $4$-secant line.
\vskip 6pt
\indent \textbf{i)} Let us settle first the case $(d',g')=(d+1,g+1)$. Take $p,q\in C$ general points, $L=\overline{pq}\subseteq \mathbb P^3$ and $X:=C\cup L$. By applying Lemma 1.2 from [BE], we know that  $H^1(X,N_X)=0$ and the curve $X$ is smoothable in $\mathbb P^3$, that is, there exists a flat family of curves $\{X_t\}$ in $\mathbb P^3$ over a smooth and irreducible base, with the general fibre $X_t$ smooth while the special fibre $X_0$ is $X$. If $d-4<[(g+3)/2]$, then since $C$ has only finitely many
$\mathfrak g^1_{d-4}$'s, by applying Prop.4.2 we get that $\mbox{gon}(X)=d-3$. In the case $d-4\geq [(g+3)/2]$ we just notice that $\mbox{gon}(X)\geq \mbox{gon}(C)=[(g'+3)/2]$.
\vskip 6pt
\indent \textbf{ii)} Next we tackle the case $(d',g')=(d+1,g+2)$. Assume first that $d-4<[(g+3)/2]\Leftrightarrow d'-4<[(g'+3)/2]$. We apply Lemma 1.2 from [BE] to a curve $X:=C\cup L$, where $L$ is a suitable trisecant line to $C$. In order to conclude that $X$ is smoothable in $\mathbb P^3$ and that $H^1(X,N_X)=0$,
we have to make sure that the trisecant line $L=\overline{pqq'}$ with $p,q,q'\in C$ can be chosen in such a way that
\begin{equation}
L, T_p(C), T_{q}(C)\mbox{ and }T_{q'}(C) \mbox{ do not all lie in the same plane. }
\end{equation}
\noindent We claim that when $C\in |C|$ is general in its linear system, at least one of its trisecants satisfies (3). Suppose not. Then for {\sl every}
smooth curve $C\in |C|$ and for {\sl every} trisecant line $L$ to $C$ condition (3) fails.
\newline
\indent We consider a $0$-dimensional subscheme $Z\subseteq S$ where $Z=p+q+q'+u+u'$, with $p,q,q'\in S$ being collinear points while $u$ and $u'$ are general infinitely near points to $q$ and $q'$ respectively. The linear system $|C|$ is at least $5$-very ample (cf. Remark 1), hence a general curve $C\in |C-Z|$ is smooth and possesses a trisecant line for which (3) holds, a contradiction. \newline
\indent Since the scheme of trisecants to a space curve is of pure dimension $1$, it follows that for a general curve $C\in |C|$, through a general point $p\in C$ there passes a trisecant line $L$
for which (3) holds. We have that $X:=C\cup L$ is smoothable in $\mathbb P^3$     and $H^1(X,N_X)=0$. We conclude that $\mbox{gon}(X)=d-3$ by proving that there is no $\mathfrak g^1_{d-4}$
on $C$ containing $L\cap C$ in a fibre.
\newline
\indent If $C\in |C|$ is general, any line in $\mathbb P^3$ (hence also a $4$-secant line to $C$) can meet only finitely many trisecants. Indeed, assuming that $m\subseteq \mathbb P^3$ is a line meeting infinitely many trisecants,
we consider the correspondence 
$$T=\{(p,t)\in C\times m:\overline{pt}\mbox{ is a trisecant to }C\}$$
and the projections $\pi_1:T\rightarrow C$ and $\pi_2:T\rightarrow m$. If $\pi_2$ is surjective, then $\mbox{Nm}_{\pi_1}(\pi_2)$ yields a $\mathfrak g^1_3$ on $C$, a contradiction. If $\pi_2$ is not surjective then there exists a point 
$t\in \mathbb P^3$ such that $\overline{pt}$ is a trisecant to $C$ for each $p\in C$. This possibility cannot occur for a general $C\in |C|$: Otherwise we take general points $t\in \mathbb P^3$ and $p,p'\in S$ and
if we denote 
$$\mathcal{B}:=\{C\in |C|:p,p'\in C\mbox{ and }\overline{tx} \mbox{ is a trisecant to }C\mbox{ for each }x\in C\},$$
we have that $\mbox{dim }\mathcal{B}\geq g-5$. On the other hand since $\overline{tp}$ and $\overline{tp'}$ are trisecants for all curves $C\in \mathcal{B}$, there must be a $0$-dimensional subscheme $Z\subseteq (\overline{tp}\cup \overline{tp'})\cap S$ of length $6$ such that $\mathcal{B}\subseteq |C-Z|$, hence $\mbox{dim}\ \mathcal{B}\leq \mbox{dim}
|C-Z|=g-6$ (use again that $|C|$ is $5$-very ample), a contradiction.
In this way the case $d-4<[(g+3)/2]$ is settled.
\newline
\indent When $d-4\geq [(g+3)/2]$ we apply Theorem 3 to obtain a smooth curve $C_1\subseteq \mathbb P^3$ of degree $d$ and genus $g+2$ such that $\mbox{gon}(C_1)=(g+5)/2$ and $H^1(C_1,N_{C_1})=0$. We take $X_1:=C_1\cup L_1$ with $L_1$ being a general $1$-secant line to $C_1$. Then $X_1$ is smoothable and $\mbox{gon}(X_1)=\mbox{gon}(C_1)=(g+5)/2$.
\vskip 6pt
\indent \textbf{iii)} Finally, we turn to the case $(d',g')=(d+2,g+3)$. Take $H\subseteq \mathbb P^3$ a general plane meeting $C$ in $d$ distinct points in general linear position and pick 4 of them: $p_1,p_2,p_3,p_4\in C\cap H$. Choose $Q\subseteq H$ a general conic such that $Q\cap C=\{p_1,p_2,p_3,p_4\}$. Theorem 5.2 from [Se] ensures that $X:=C\cup Q$ is smoothable in $\mathbb P^3$ and $H^1(X,N_X)=0$.\newline 
\indent Assume first that $d'-4\leq [(g'+3)/2]$. We claim that $\mbox{gon}(X)\geq \mbox{gon}(C)+2$. According to Prop.4.2 the opposite could happen only in 2 cases: a) There exists a $\mathfrak g^1_{d-3}$ on $C$, say $|Z|$,
such that $|Z|(-p_1-p_2-p_3-p_4)\neq \emptyset.$\  \ b) There exists a degree $d-4$ map $f:C\rightarrow \mathbb P^1$ and a degree $1$ map $f':Q\rightarrow \mathbb P^1$ such that $f(p_i)=f'(p_i)$, for $i=1,\ldots, 4$.
\newline
\indent Assume that a) does happen. We denote by $U=\{D\in C_4:|\mathcal{O}_C(1)|(-D)\neq \emptyset\}$ the irreducible $3$-fold of
divisors of degree $4$ spanning a plane and also consider the correspondence 
$$\Sigma=\{(L,D)\in W^1_{d-3}(C)\times U:|L|(-D)\neq \emptyset\},$$
with the projections $\pi_1:\Sigma\rightarrow W^1_{d-3}(C)$ and $\pi_2:\Sigma\rightarrow U$. We know that $\pi_2$ is dominant, hence $\mbox{dim }\Sigma\geq 3$ and therefore $\mbox{dim }W^1_{d-3}(C)\geq 2$. 
\newline
\indent If $\rho(g,1,d-3)<0$ by Prop.3.1 we get that every base-point-free $\mathfrak g^1_{d-3}$ on $C$ is cut out by a divisor $D$ on $S$ such that $D\in \mathcal{A}$ (see the proof of Theorem 3 for this notation) and $C\cdot D-D^2=\mbox{Cliff}(C,D_{|C})+2\leq d-3$, hence $C\cdot D-D^2\leq d-4$ for parity reasons. As pointed out at the end of the proof of Theorem 3 this forces $D\sim H$, that is, all base-point-free $\mathfrak g^1_{d-3}$'s on $C$ are given by planes through a trisecant line.  Thus $C$ has  $\infty^2$ trisecants, a contradiction.  
\newline \indent If $\rho(g,1,d-3)\geq 0$, then $g=2d-9$ and we can assume that there is $L\in \pi_1(\Sigma)$ such that $|\mathcal{O}_C(1)-L|\neq \emptyset$. The map $\pi_1$ is either generically finite hence $\mbox{dim }W^1_{d-4}(C)\geq \mbox{dim }W^1_{d-3}(C)-2\geq 1$ (cf. [FHL]), a contradiction, or otherwise $\pi_1$ has fibre dimension $1$. This is possible only when there is a component $A$ of $W^1_{d-3}(C)$ with $\mbox{dim}(A)\geq 2$ and such that the general $L\in A$ satisfies $|\mathcal{O}_C(1)-L|\neq \emptyset$ and {\sl every} $L\in A$ has non-ordinary ramification so that the monodromy of each $\mathfrak g^1_{d-3}$ is not the full symmetric group. Applying again [FHL] there is $L\in W^1_{d-4}(C)$ such that $\{L\}+W^0_1(C)\subseteq A$, in particular $L$ has non-ordinary ramification too. It is easy to see that this contradicts the $(d-7)$-very ampleness of $|C|$ asserted by Remark 1.
\newline
\indent We now rule out case b). Suppose that b) does happen and denote by
$L\subseteq \mathbb P^3$ the $4$-secant line corresponding to $f$. Let $\{p\}=
L\cap H$, and pick $l\subseteq H$ a general line. As $Q$ was a general conic
through $p_1,\ldots ,p_4$ we may assume that $p\notin Q$. The map $f':Q\rightarrow l$ is (up to a projective isomorphism of $l$) the projection from a point $q\in Q$, while $f(p_i)=\overline{p_ip}\cap l$, for $i=1,\ldots,4$.
By Steiner's Theorem from classical projective geometry, the condition 
$(f(p_1)f(p_2)f(p_3)f(p_4))=(f'(p_1)f'(p_2)f'(p_3)f'(p_4))$ is equivalent with 
$p_1,p_2,p_3,p_4, p$ and $q$ being on a conic, a contradiction since $p\notin Q$.\newline
\indent Finally, when $d'-4>[(g'+3)/2]$, we have to show that $\mbox{gon}(X)\geq \mbox{gon}(C)+1$. We note that $\mbox{dim }G^1_{(g+3)/2}(C)=1$ (for any curve
one has the inequality
$\mbox{dim }G^1_{\tiny{\mbox{gon}}}\leq 1$). By taking $H\in (\mathbb P^3)^{\vee}$ general enough, we obtain that $p_1,\ldots ,p_4$ do not occur in the same fibre of a $\mathfrak g^1_{(g+3)/2}$.\hfill $\Box$
\vskip 4pt
\noindent \textbf{Remark: } Theorem 1 can be viewed as a non-containment relation
$\mathcal{M}_{g',d'}^3\nsubseteq \mathcal{M}_{g',d'-5}^1$ between different Brill-Noether loci when $d'$ and $g'$ are as in Theorem 1 and moreover $d'-4\leq [(g'+3)/2]$. We can turn this problem on its head and ask the following question: given $g$ and $k$ such that $k<(g+2)/2$, when is it true that the general $k$-gonal curve of genus $g$ has no other linear series $\mathfrak g^r_d$ with $\rho(g,r,d)<0$, that is, the pencil computing the gonality is the only Brill-Noether exceptional linear series?
\newline
\indent In [Fa2] we prove using limit linear series the following result:
fix $g$ and $k$ positive integers such that $-3\leq \rho(g,1,k)<0$. If $\rho(g,1,k)=-3$ assume furthermore that $k\geq 6$. Then the general $k$-gonal curve $C$ of genus $g$ has no $\mathfrak g^r_d$'s with $\rho(g,r,d)<0$ except
$\mathfrak g^1_k$ and $|K_C-\mathfrak g^1_k|$. In other words the $k$-gonal locus $\mathcal{M}_{g,k}^1$ is not contained in any other proper Brill-Noether locus $\mathcal{M}_{g,d}^r$ with $r\geq 2,d\leq g-1$ and $\rho(g,r,d)<0$.
\newline
\indent In seems that other methods are needed to extend this result for more negative values of $\rho(g,1,k)$.

\section{The Kodaira dimension of $\mathcal{M}_{23}$}
In this section we explain how Theorem 1 gives a new proof of our result $\kappa(\mathcal{M}_{23})\geq 2$ (cf. [Fa]). We refer to [Fa] for a detailed analysis of the geometry of $\mathcal{M}_{23}$ ; in that paper we also conjecture that $\kappa(\mathcal{M}_{23})=2$ and we present evidence for such a possibility.
\newline
\indent Let us denote by $\overline{\mathcal{M}}_{g}$ the moduli space of Deligne-Mumford stable curves of genus $g$. We study the multicanonical linear systems on $\overline{\mathcal{M}}_{23}$ by exhibiting three explicit multicanonical divisors on $\overline{\mathcal{M}}_{23}$ which are (modulo a
positive combination of boundary classes coming from $\overline{\mathcal{M}}_{23}-\mathcal{M}_{23}$) of Brill-Noether type, that is,
loci of curves having a $\mathfrak g^r_d$ when $\rho(23,r,d)=-1$.
\newline
\indent On $\mathcal{M}_{23}$ there are three Brill-Noether divisors corresponding to the solutions of the equation $\rho(23,r,d)=-1$: the $12$-gonal divisor $\mathcal{M}^1_{23,12}$, the divisor $\mathcal{M}^2_{23,17}$ of curves having a $\mathfrak g^2_{17}$ and finally the divisor $\mathcal{M}_{23,20}^3$ of curves possessing a $\mathfrak g^3_{20}$. If we denote by $\overline{\mathcal{M}}_{g,d}^r$ the closure of $\mathcal{M}_{g,d}^r$ inside $\overline{\mathcal{M}}_g$, the classes $[\overline{\mathcal{M}}_{g,d}^r]\in \mbox{Pic}_{\mathbb Q}(\overline{\mathcal{M}}_g)$ when $\rho(g,r,d)=-1$ have been computed (see [EH],[Fa]). It is quite remarkable that for fixed $g$ all classes $[\overline{\mathcal{M}}_{g,d}^r]$ are proportional. One also knows the canonical divisor class (cf. [HM]):
$$K_{\overline{\mathcal{M}}_g}=13\lambda-2\delta_0-3\delta_1-2\delta_2-\cdots -2\delta_{[g/2]},$$
and by comparing for $g=23$ this formula with the expression of the classes $[\overline{\mathcal{M}}_{23,d}^r]$, we find that there are constants
$m,m_1,m_2,m_3\in \mathbb Z_{>0}$ such that
$$mK_{\overline{\mathcal{M}}_{23}}=m_1[\overline{\mathcal{M}}^1_{23,12}]+E=m_2[\overline{\mathcal{M}}_{23,17}^2]+E=m_3[\overline{\mathcal{M}}_{23,20}^3]+E,$$
where $E$ is the same positive combination of the boundary classes $\delta_1,\ldots ,\delta_{11}$.
\newline
\indent As explained in [Fa], since    $\overline{\mathcal{M}}_{23,12}^1$,\mbox{ }$\overline{\mathcal{M}}_{23,17}^2$ and 
$\overline{\mathcal{M}}_{23,20}^3$ are mutually distinct 
irreducible divisors, we can show that the multicanonical image of $\overline{\mathcal{M}}_{23}$ cannot be a curve once 
we construct a smooth curve of genus $23$ lying in the 
support of exactly two of the divisors $\mathcal{M}_{23,12}^1$,$\mathcal{M}^2_{23,17}$ and $\mathcal{M}^3_{23,20}$.
\ In this way we rule out the possibility of all three intersections of two Brill-Noether divisors being equal to 
base-locus$(|mK_{\overline{\mathcal{M}}_{23}}|)\cap \mathcal{M}_{23}$.
\newline
\indent In [Fa] we found such genus $23$ curves using 
an intricate construction involving limit linear series 
(cf. Proposition 5.4 in [Fa]). Here we can construct 
such curves in a much simpler way by applying Theorem 1 
when $(d,g)=(18,23)$: there exists a smooth curve 
$C\subseteq \mathbb P^3$ of genus $23$ and degree $18$ such that $\mbox{gon}(C)=13$; moreover $C$ sits on a smooth quartic surface $S\subseteq \mathbb P^3$ such that $\mbox{Pic}(S)=\mathbb ZC\oplus \mathbb Z H.$
\newline
\indent Since $C$ has a very ample $\mathfrak g^3_{18}$, by adding $2$ base points it will also have plenty of $\mathfrak g^3_{20}$'s and also $\mathfrak g^2_{17}$'s of the form $\mathfrak g^3_{18}(-p)=\{D\in \mathfrak g^3_{18}:D\geq p\}$, for any $p\in C$. Therefore $[C]\in (\mathcal{M}_{23,20}^3\cap \mathcal{M}_{23,17}^2)-\mathcal{M}_{23,12}^1$, and Theorem 2 now follows.

\vskip10pt
\noindent Department of Mathematics, University of Michigan\newline
525 East University, Ann Arbor, MI 48109
\newline
e-mail: {\tt gfarkas@math.lsa.umich.edu }


{\footnotesize
\begin{thebibliography}{[AC]}
\bibitem[AC]{AC} E. Arbarello, \ M. Cornalba, \ {\em{Footnotes  to a paper of
Beniamino Segre, }} Math. Ann. 256(1981), 341--362.
\bibitem[ACGH]{ACGH} E.\ Arbarello, \ M.\ Cornalba, \ P.\ Griffiths, \ J.\ Harris,  \
{\em{Geometry of algebraic curves, }} Grundlehren der Mathematischen Wissenschaften,
267, Springer Verlag, 1985.
\bibitem[BE]{BE} E. Ballico, \ Ph. Ellia, \ {\em{On the existence of curves with maximal rank in $\mathbb P^n$, }} J.\ reine angew. Math. 397(1989), 1-22.
\bibitem [Ba]{Ba} B. Basili, \ {\em{Indice de Clifford des intersections completes de l'espace, }} Bull. Soc. Math. France 124(1996), no.1, 61-95.
\bibitem[BS]{BS} M. Beltrametti, \ A. Sommese, \ {\em{Zero cycles and $k$-th order embeddings of smooth projective surfaces, }} in: Problems in the theory of surfaces and their classification (Cortona, 1988), 33-48, Academic Press, London, 1991.
\bibitem[CilP]{CilP} C. Ciliberto, \ G. Pareschi, \ {\em{Pencils of minimal degree on curves on a $K3$ surface, }} J. reine angew. Math., 460(1995), 15-36.
\bibitem [Co]{Co} M. Coppens, \ {\em{The gonality of general smooth curves with a prescribed plane model, }} Math. Ann. 289(1991), 89-93.
\bibitem [DMo]{DMo} R. Donagi, \ D. Morrison, \ {\em{Linear systems on $K3$ sections, }} J. Diff. Geom. 29(1989), 49-64.
\bibitem[EH]{EH} D.\ Eisenbud, \ J. Harris, \ {\em{The Kodaira dimension  of the
moduli space of curves of genus $\geq 23,$ }} Invent.\ Math.\ 90(1987), no.2, 
359--387.
\bibitem[Fa]{Fa} G.\ Farkas, \ {\em{The geometry of the moduli space of curves of genus $23$, }} Math. Ann. 308(2000), no.1, 43-65.
\bibitem[Fa2]{Fa2} G.\ Farkas, \ {\em{The birational geometry of the moduli space of curves, }} Ph. D. Thesis, Universiteit van Amsterdam, 2000.
\bibitem[FHL]{FHL} W. Fulton, \ J. Harris, \ R. Lazarsfeld, \ {\em{Excess linear series on an algebraic curve, }} Proc. AMS 92(1984), 320-322.
\bibitem[GH]{GH} P. Griffiths, \ J. Harris, \ {\em{Infinitesimal variations of Hodge structure: an infinitesimal invariant of Hodge classes, }} Compositio Math. 50(1983), 207-265.
\bibitem [GL]{GL} M. Green, \ R. Lazarsfeld, \ {\em{Special divisors on curves on $K3$ surfaces, }} Invent. Math. 89(1987), 73-90.
\bibitem [HM]{HM} J.\ Harris, \ D.\ Mumford, \ {\em{On the Kodaira dimension of  the
moduli space of curves,}} Invent. Math. 67(1982), no.1, 23--88.
\bibitem [Mod]{Mod} J. Harris, \ I. Morrison, \ {\em{Moduli of curves, }}  Graduate
Texts in Mathematics, 187, Springer Verlag, 1998.
\bibitem [Hir]{Hir} A. Hirschowitz, \ {\em{Section planes et multisecantes pour les courbes gauches generiques principales, }} in: Space curves: Proceedings, Rocca di Pappa 1985,
\ Lecture Notes in Math. 1226, Springer Verlag 1987, 124-155.
\bibitem [Kn]{Kn}  A. Knutsen, \ {\em{On degrees and genera of smooth curves on projective $K3$ surfaces, }} math.AG/98050140 preprint.
\bibitem [Laz]{Laz} R. Lazarsfeld, \ {\em{Brill-Noether-Petri without degenerations, }} J. Diff. Geom. 23(1986), 299-307.
\bibitem [Mo]{Mo} S. Mori, \ {\em{On degrees and genera of curves on smooth quartic surfaces in $\mathbb P^3$, }} Nagoya Math. J. 96(1984), 127-132.
\bibitem [Ra]{Ra} J. Rathmann, \ {\em{The genus of algebraic space curves, }} Ph. D. Thesis, Berkeley, 1986.
\bibitem[SD]{SD} B.\ Saint-Donat, \ {\em{Projective models of $K3$ surfaces, }}
Amer. J. Math., 96(1974), 602-639.
\bibitem [Se]{Se} E.\ Sernesi, \ {\em{On the existence of certain families of 
curves, }} Invent.\ Math.\ 75(1984), no.1, 25--57.
\end{thebibliography}}
\end{document}